\pdfoutput=1
\PassOptionsToPackage{pdfusetitle,psdextra}{hyperref}
\PassOptionsToPackage{style=numeric-comp}{biblatex}
\PassOptionsToPackage{capitalize}{cleveref}
\documentclass[a4paper]{article}

\usepackage{xr}
\usepackage{subfiles}
\usepackage{biblatex}
\usepackage{amsthm}
\usepackage{my-preamble}
\usepackage{xassoccnt}

\NewCommandCopy{\myCite}{\parencite}
\NewCommandCopy{\myTextCite}{\textcite}

\makeatletter
\NewDocumentCommand{\counter@within}{}{section}
\DeclareCoupledCountersGroup{thm@group}
\NewDocumentCommand{\myNewTheorem}{m m}{%
  \newtheorem{#1}{#2}[\counter@within]
  \AddCoupledCounters[name=thm@group]{#1}
  \Crefname{#1}{#2}{#2s}
}
\makeatother
\theoremstyle{definition}
\myNewTheorem{myDefinition}{Definition}
\myNewTheorem{myConstruction}{Construction}
\myNewTheorem{myNotation}{Notation}
\theoremstyle{remark}
\myNewTheorem{myExample}{Example}
\myNewTheorem{myRemark}{Remark}
\theoremstyle{plain}
\myNewTheorem{myProposition}{Proposition}
\myNewTheorem{myCorollary}{Corollary}
\myNewTheorem{myLemma}{Lemma}
\myNewTheorem{myTheorem}{Theorem}
\myNewTheorem{myFact}{Fact}

\NewDocumentEnvironment{myProof}{}{\begin{proof}}{\end{proof}}
\NewCommandCopy{\myQed}{\qed}

\DeclareNameAlias{author}{given-family}
\DeclareFieldFormat{postnote}{#1}

\author{\myAuthor}
\title{\myTitle}
\date{June 12, 2025}

\begin{document}

\maketitle

\begin{abstract}

We provide an alternative proof of Lurie's result that the wide subcategory of the $\infty$-category of $\infty$-topoi spanned by the étale morphisms is closed under small colimits.
Our proof is based on a new characterization of étale morphisms of $\infty$-topoi in relation to univalent families and does not rely on a larger universe.
During the proof, we also give an elementary construction of univalent completion.

\end{abstract}

\section{Introduction}
\label{sec:introduction}

An \emph{\(\infty\)-topos} is a ``generalized space'' whose behavior is completely determined by the \(\infty\)-category of sheaves on it.
\(\infty\)-topoi are an ideal place for homotopy theory, homology theory, and cohomology theory.
\emph{Étale morphisms} play the central role in \(\infty\)-topos theory in that the \(\infty\)-category \(\mySh(X)\) of sheaves on an \(\infty\)-topos \(X\) is equivalent to the \(\infty\)-category of étale morphisms into \(X\).

In his book, Lurie showed that the wide subcategory of the \(\infty\)-category of \(\infty\)-topoi spanned by the étale morphisms is closed under small colimits \myCite[Theorem 6.3.5.13]{lurie2009higher}.
His proof is, however, rather technical and uses one larger universe.
The purpose of this paper is to provide an alternative proof of Lurie's result without assuming a larger universe.
Our proof is based on a new characterization of étale morphisms in relation to \emph{univalent families} (\cref{prop-bc-char-1}).

There is a known characterization of étale morphisms as those morphisms \(f \myElemOf Y \myMorphism X\) of \(\infty\)-topoi whose inverse image functor \(\myInvIm{f} \myElemOf \mySh(X) \myMorphism \mySh(Y)\) has a left adjoint that satisfies certain conditions \myCite[Proposition 6.3.5.11]{lurie2009higher}.
However, this characterization is not so helpful in proving the closure property of étale morphisms because the conditions for the left adjoint of \(\myInvIm{f}\) are hard to verify unless the left adjoint is explicitly calculated.

A family (i.e. a morphism) in \(\mySh(X)\) is called univalent if there is at most one pullback to it from any family \myCite[Proposition 3.8]{gepner2017univalence}.
It is a characteristic property of the \(\infty\)-categories of sheaves on \(\infty\)-topoi that, for every family \(u\), there exists a univalent family \(v\) such that \(u\) is a pullback of \(v\) \myCite[Theorem 6.1.6.8]{lurie2009higher}.

We characterize étale morphisms by properties of \(\myInvIm{f}\) itself (\cref{prop-bc-char-1}):
A morphism \(f \myElemOf Y \myMorphism X\) of \(\infty\)-topoi is étale if and only if its inverse image functor \(\myInvIm{f} \myElemOf \mySh(X) \myMorphism \mySh(Y)\) preserves dependent products (i.e. right adjoints of pullback functors) and satisfies that, for any univalent family \(v\) in \(\mySh(Y)\), there exists a univalent family \(u\) in \(\mySh(X)\) such that \(v\) is a pullback of \(\myInvIm{f}(u)\).

\subsection{Organization}

In \cref{sec:category-theory}, we collect some facts in general \(\infty\)-category theory.
We review basic \(\infty\)-topos theory in \cref{sec:topoi-logoi}.
In \cref{sec:univalence-logoi}, we study univalent families in \(\infty\)-categories of sheaves on \(\infty\)-topoi.
We provide an elementary construction of the \emph{univalent completion} of an arbitrary family (\cref{prop-univalent-completion}), that is, the least univalent family that extends the given family.
Finally in \cref{sec:etale-morph-topoi}, we prove our characterization of étale morphisms of \(\infty\)-topoi (\cref{prop-bc-char-1}) and give an alternative proof of Lurie's result that the wide subcategory of the \(\infty\)-category of \(\infty\)-topoi spanned by the étale morphisms is closed under small colimits (\cref{prop-etale-colimit}).

\subsection{Foundations}

Our formal metatheory is {\rm ZFC} with one Grothendieck universe of small sets, but we mainly work in the homotopy theory of Kan complexes.
Every construction is thus homotopy invariant.
We refer to Kan complexes as \myDefine{types} and prefer type-theoretic notation \(a \myElemOf A\) to set-theoretic notation \(a \in A\) to mean that \(a\) is an element of a type \(A\).
To put it more formally, one could use type theory with the univalence axiom \myCite{hottbook} and interpret it in Kan complexes \myCite{kapulkin2021simplicial}.
The only reason for the non-constructive metatheory is that existing theories of \(\infty\)-categories and presentable \(\infty\)-categories are non-constructive.
We drop the prefix ``\(\infty\)-'', so \(\infty\)-categories and \(\infty\)-topoi are simply called categories and topoi, respectively.

\section{Category theory}
\label{sec:category-theory}

We collects some facts in general category theory.

Many structures on categories are defined by adjunctions.
Preservation of such structures is defined by the so-called Beck-Chevalley condition.

\begin{myDefinition}
  Consider a commutative square of categories
  \begin{equation*}
    \begin{tikzcd}
      D
      \arrow[r, "G"]
      \arrow[d, "P"'] &
      D'
      \arrow[d, "P'"] \\
      C
      \arrow[r, "F"'] &
      C'.
    \end{tikzcd}
  \end{equation*}
  We view \((G \myComma F)\) as a morphism (in the category of categories) from \(P\) to \(P'\), leaving the homotopy filling the square implicit.
  Let \(\myAdjUnit\) and \(\myAdjCounit\) denote the unit and counit, respectively, of an adjunction.
  \begin{itemize}
  \item Suppose that \(P\) and \(P'\) have left adjoints \(Q\) and \(Q'\), respectively.
    We say the morphism \((G \myComma F)\) from \(P\) to \(P'\) satisfies the \myDefine{left Beck-Chevalley condition} if the natural transformation \(Q' \myComp F \myMMorphism G \myComp Q\) defined by the composite
    \begin{equation*}
      \begin{tikzcd}
        & D
        \arrow[rr, "G"]
        \arrow[dr, "P"] & &
        D'
        \arrow[rr, equal]
        \arrow[dr, "P'"'] &
        \phantom{D'} &
        D' \\
        C
        \arrow[ur, "Q"]
        \arrow[rr, equal] &
        \phantom{C}
        \arrow[u, myMMorphism, "\myAdjUnit"] &
        C
        \arrow[rr, "F"'] & &
        C'
        \arrow[ur, "Q'"']
        \arrow[u, myMMorphism, "\myAdjCounit"']
      \end{tikzcd}
    \end{equation*}
    is an equivalence.
  \item Suppose that \(P\) and \(P'\) have right adjoints \(R\) and \(R'\), respectively.
    We say the morphism \((G \myComma F)\) from \(P\) to \(P'\) satisfies the \myDefine{right Beck-Chevalley condition} if the natural transformation \(G \myComp R \myMMorphism R' \myComp F\) defined by the composite
    \begin{equation*}
      \begin{tikzcd}
        & D
        \arrow[rr, "G"]
        \arrow[dr, "P"]
        \arrow[d, myMMorphism, "\myAdjCounit"'] & &
        D'
        \arrow[rr, equal]
        \arrow[dr, "P'"'] &
        \phantom{D'}
        \arrow[d, myMMorphism, "\myAdjUnit"] &
        D' \\
        C
        \arrow[ur, "R"]
        \arrow[rr, equal] &
        \phantom{C} &
        C
        \arrow[rr, "F"'] & &
        C'
        \arrow[ur, "R'"']
      \end{tikzcd}
    \end{equation*}
    is an equivalence.
  \end{itemize}
\end{myDefinition}

\begin{myNotation}
  Let \(C\) be a category and \(f \myElemOf x \myMorphism x'\) a morphism in \(C\).
  We refer to the functor \((y \myComma p) \myMapsTo (y \myComma f \myComp p) \myElemOf C \mySlice x \myMorphism C \mySlice x'\) as \(\mySumFun{f}\).
  When \(C\) has pullbacks along \(f\), let \(\myPBFun{f} \myElemOf C \mySlice x' \myMorphism C \mySlice x\) denote the pullback functor.
  Note that \(\mySumFun{f} \myAdjRel \myPBFun{f}\).
  When \(\myPBFun{f}\) further has a right adjoint, we refer to it as \(\myProdFun{f} \myElemOf C \mySlice x \myMorphism C \mySlice x'\) and call \(\myProdFun{f}(y)\) the \myDefine{dependent product of \(y\) along \(f\)} for \(y \myElemOf C \mySlice x\).
  When \(x'\) is the terminal object in \(C\) and \(f\) is the unique morphism to it, we write \(\mySumFun{x}\), \(\myPBFun{x}\), and \(\myProdFun{x}\) for \(\mySumFun{f}\), \(\myPBFun{f}\), and \(\myProdFun{f}\), respectively.
\end{myNotation}

\begin{myExample}
  Let \(F \myElemOf C \myMorphism D\) be a functor between categories with pullbacks.
  For an object \(x \myElemOf C\), let \(F \mySlice x \myElemOf C \mySlice x \myMorphism D \mySlice F(x)\) denote the functor defined by \((F \mySlice x)(y \myComma p) \myDefEq (F(y) \myComma F(p))\).
  Then, \(F\) preserves pullbacks if and only if, for every morphism \(f \myElemOf x \myMorphism x'\) in \(C\), the morphism \((F \mySlice x \myComma F \mySlice x')\) from \(\mySumFun{f}\) to \(\mySumFun{F(f)}\) satisfies the right Beck-Chevalley condition.
\end{myExample}

\begin{myDefinition}
  Let \(F \myElemOf C \myMorphism D\) be a functor between categories with pullbacks and dependent products that preserves pullbacks.
  We say \(F\) \myDefine{preserves dependent products} if, for every morphism \(f \myElemOf x' \myMorphism x\) in \(C\), the morphism \((F \mySlice x \myComma F \mySlice x')\) from \(\myPBFun{f}\) to \(\myPBFun{F(f)}\) satisfies the right Beck-Chevalley condition.
\end{myDefinition}

The left and right Beck-Chevalley conditions are related to each other.

\begin{myProposition}
  \label{prop-left-right-bc}
  Let
  \begin{equation*}
    \begin{tikzcd}
      A
      \arrow[r, "F"]
      \arrow[d, "G"'] &
      B
      \arrow[d, "H"] \\
      C
      \arrow[r, "K"'] &
      D
    \end{tikzcd}
  \end{equation*}
  be a commutative square of categories.
  Suppose that \(F\) and \(K\) have left adjoints and \(G\) and \(H\) have right adjoints.
  Then the morphism \((G \myComma H)\) from \(F\) to \(K\) satisfies the left Beck-Chevalley condition if and only if the morphism \((F \myComma K)\) from \(G\) to \(H\) satisfies the right Beck-Chevalley condition.
\end{myProposition}
\begin{myProof}
  Let \(F'\) and \(K'\) be the left adjoints of \(F\) and \(K\), respectively, and \(G'\) and \(H'\) the right adjoints of \(G\) and \(H\), respectively.
  Let \(X \mySub C^{B}\) denote the full subcategory spanned by the left adjoint functors, and \(Y \mySub B^{C}\) the full subcategory spanned by the right adjoint functors.
  By definition, \(\myOpCat{X} \myEquiv Y\).
  The natural transformation \(K' \myComp H \myMMorphism G \myComp F'\) corresponds to the natural transformation \(F \myComp G' \myMMorphism H' \myComp K\) via this equivalence.
  Hence, one is an equivalence if and only if the other is.
\end{myProof}

We will use the following criterion for a wide subcategory to be closed under colimits.

\begin{myProposition}
  \label{prop-closed-under-colimit-criterion-1}
  Let \(C\) be a category with small colimits and \(D \mySub C\) a wide subcategory.
  Suppose:
  \begin{enumerate}
    \NewCommandCopy{\myLabelOrig}{\label}
    \RenewDocumentCommand{\label}{m}{\myLabelOrig[cond]{#1}}
  \item \label{cond-closed-under-colimit-criterion-1-slice}
    \(D \mySlice x \mySub C \mySlice x\) is closed under small colimits for every \(x \myElemOf C\).
  \item \label{cond-closed-under-colimit-criterion-1-proj}
    For every small diagram \(x \myElemOf I \myMorphism D\) with colimit \(x_{\infty}\) in \(C\), all the inclusions \(\myColimInj{i} \myElemOf x_{i} \myMorphism x_{\infty}\) are in \(D\).
  \end{enumerate}
  Then \(D \mySub C\) is closed under small colimits.
\end{myProposition}
\begin{myProof}
  Let \(x \myElemOf I \myMorphism D\) be a small diagram.
  Take the colimit \(x_{\infty}\) of \(x\) in \(C\).
  It suffices to show that a morphism \(f \myElemOf x_{\infty} \myMorphism y\) belongs to \(D\) if and only if all the inclusions \(f \myComp \myColimInj{i} \myElemOf x_{i} \myMorphism y\) belong to \(D\).
  The ``if'' direction follows from \cref{cond-closed-under-colimit-criterion-1-slice} since \(f\) is the colimit in \(C \mySlice y\) of \(i \myMapsTo f \myComp \myColimInj{i}\).
  The ``only if'' direction follows from \cref{cond-closed-under-colimit-criterion-1-proj}.
\end{myProof}

\section{Topoi and logoi}
\label{sec:topoi-logoi}

We recall the notion of a topos \myCite[Chapter 6]{lurie2009higher}.
Since the goal of this paper is an alternative proof of a theorem in topos theory, we do not rely on any result in topos theory to avoid a circular argument.

We assume that the reader is familiar with the theory of presentable categories \myCite[Section 5.5]{lurie2009higher}.
A category \(C\) is \myDefine{presentable} if it has small colimits and there exist a small category \(A\) and a small regular cardinal \(\kappa\) such that \(C\) is equivalent to a certain category \(\myInd_{\kappa}(A)\).
Take a cardinal \(\lambda\) greater than the cardinalities of \(\myInd_{\kappa}(A)\) for all small categories \(A\) and small regular cardinals \(\kappa\), and then we define \(\myPrL\) to be the subcategory of the category of \(\lambda\)-small categories spanned by the presentable categories and the functors preserving small colimits.
We use the following facts on presentable categories.

\begin{myFact}[{\myTextCite[Proposition 5.5.3.13]{lurie2009higher}}]
  \label{prop-pres-l-limit}
  \(\myPrL\) has small limits computed as limits of categories.
\end{myFact}

\begin{myFact}[{\myTextCite[Proposition 5.5.2.2]{lurie2009higher}}]
  \label{prop-pres-left-repr}
  A presheaf of small types on a presentable category is representable if and only if it preserves small limits.
\end{myFact}

\begin{myFact}[{\myTextCite[Corollary 5.5.2.9 (1)]{lurie2009higher}}]
  \label{prop-pres-left-aft}
  A functor between presentable categories is a left adjoint if and only if it preserves small colimits.
\end{myFact}

\begin{myFact}[{\myTextCite[Corollary 5.5.2.9 (2)]{lurie2009higher}}]
  \label{prop-pres-right-aft}
  A functor \(F\) between presentable categories is a right adjoint if and only if \(F\) preserves small limits and there exists a small regular cardinal \(\kappa\) such that \(F\) preserves small \(\kappa\)-filtered colimits.
\end{myFact}

\begin{myFact}[{\myTextCite[Proposition 5.5.3.10]{lurie2009higher}}]
  \label{prop-pres-slice}
  For every presentable category \(C\) and object \(x \myElemOf C\), the slice \(C \mySlice x\) is presentable.
\end{myFact}

We follow \myTextCite{anel2021topo-logie} and use the term logos for the algebraic dual of a topos.

\begin{myDefinition}
  We say a natural transformation is \myDefine{cartesian} if all the naturality squares are pullbacks.
\end{myDefinition}

\begin{myDefinition}
  \label{def-logos}
  A \myDefine{logos} is a presentable category \(C\) with finite limits\footnote{Finite limits are actually redundant because every presentable category has small limits \myCite[Corollary 5.5.2.4]{lurie2009higher}.} that satisfies the \myDefine{descent property}:
  For any small category \(I\) and natural transformation \(p \myElemOf y \myMMorphism x \myElemOf \myRCone{I} \myMorphism C\), if \(x\) is colimiting and \(p \myRestrict{I}\) is cartesian, then \(y\) is colimiting if and only if \(p\) is cartesian, where \(\myRCone{I}\) is the category obtained from \(I\) by freely adjoining a terminal object.
  A \myDefine{morphism of logoi} is a functor preserving small colimits and finite limits.
  Let \(\myLogos \mySub \myPrL\) denote the subcategory spanned by the logoi and the morphisms of logoi.
  The category \(\myTopos\) of \myDefine{topoi} is the opposite of \(\myLogos\).
\end{myDefinition}

We are interested in colimits of topoi (\cref{prop-colimits-of-topoi}), which are limits of logoi.

\begin{myProposition}[{\myTextCite[Proposition 6.3.2.3]{lurie2009higher}}]
  \label{prop-limits-of-logoi}
  \(\myLogos\) has small limits computed as limits of categories.
\end{myProposition}
\begin{myProof}
  Let \(C \myElemOf I \myMorphism \myLogos\) be a small diagram.
  Take the limit \(C_{-\infty}\) of categories \(C\).
  Small colimits and finite limits in \(C_{-\infty}\) exist and are computed pointwise.
  The descent property for \(C_{-\infty}\) thus follows from the descent properties for \(C_{i}\)'s.
  \(C_{-\infty}\) is presentable by \cref{prop-pres-l-limit}.
  Therefore, \(C_{-\infty}\) is the limit of \(C\) in \(\myLogos\).
\end{myProof}

\begin{myCorollary}
  \label{prop-colimits-of-topoi}
  \(\myTopos\) has small colimits.
\end{myCorollary}
\begin{myProof}
  By \cref{prop-limits-of-logoi}.
\end{myProof}

Logoi have richer structures than just small colimits and finite limits.
Dependent products are one example (\cref{prop-logos-dep-prod}).

\begin{myLemma}
  \label{prop-logos-colimit-pb-stable}
  Let \(C\) be a logos and \(f \myElemOf x' \myMorphism x\) a morphism in \(C\).
  Then the pullback functor \(\myPBFun{f} \myElemOf C \mySlice x \myMorphism C \mySlice x'\) preserves small colimits.
\end{myLemma}
\begin{myProof}
  This is a consequence of the ``if'' direction of the descent property.
\end{myProof}

\begin{myProposition}
  \label{prop-logos-dep-prod}
  Every logos has dependent products.
\end{myProposition}
\begin{myProof}
  By \cref{prop-logos-colimit-pb-stable,prop-pres-left-aft}.
\end{myProof}

\begin{myNotation}
  Let \(\myLogosLCC \mySub \myLogos\) denote the wide subcategory spanned by those morphisms that preserve dependent products.
\end{myNotation}

\begin{myProposition}
  \label{prop-logos-lcc-limit}
  \(\myLogosLCC \mySub \myLogos\) is closed under small limits.
\end{myProposition}
\begin{myProof}
  This is because dependent products in the limit in \(\myLogos\) of a small diagram \(I \myMorphism \myLogosLCC\) are computed pointwise.
\end{myProof}

The following characterization of morphisms in \(\myLogosLCC\) (\cref{prop-prsv-dp-equiv-idx-la}) is useful.

\begin{myDefinition}
  Let \(F \myElemOf C \myMorphism D\) be a morphism of logoi.
  Suppose that \(F\) has a left adjoint \(G\).
  We say \(G\) is \myDefine{\(C\)-indexed} if, for any morphism \(p \myElemOf x' \myMorphism x\) in \(C\), the morphism \((\myPBFun{p} \myComma \myPBFun{F(p)})\) from \(F \mySlice x\) to \(F \mySlice x'\) satisfies the left Beck-Chevalley condition.
  This is equivalent to that, for any pullback square in \(D\) of the form
  \begin{equation*}
    \begin{tikzcd}
      y'
      \arrow[r]
      \arrow[d, "q"']
      \arrow[dr, myPBMark] &
      F(x')
      \arrow[d, "F(p)"] \\
      y
      \arrow[r] &
      F(x),
    \end{tikzcd}
  \end{equation*}
  its transpose
  \begin{equation*}
    \begin{tikzcd}
      G(y')
      \arrow[r]
      \arrow[d, "G(q)"'] &
      x'
      \arrow[d, "p"] \\
      G(y)
      \arrow[r] &
      x
    \end{tikzcd}
  \end{equation*}
  is a pullback.
\end{myDefinition}

\begin{myLemma}[{\myTextCite[Theorem 6.1.3.9]{lurie2009higher}}]
  \label{prop-logos-self-index-sheaf}
  Let \(C\) be a logos and \(x \myElemOf I \myMorphism C\) a small diagram.
  Then the functor
  \begin{equation*}
    C \mySlice \myColim_{i \myElemOf I} x_{i} \myMorphism \myLim_{i \myElemOf I} C \mySlice x_{i}
  \end{equation*}
  defined by the pullbacks along the inclusions \(x_{i} \myMorphism \myColim_{j \myElemOf I} x_{j}\) is an equivalence.
\end{myLemma}
\begin{myProof}
  The functor
  \begin{math}
    C \mySlice \myColim_{i \myElemOf I} x_{i} \myMorphism \myLim_{i \myElemOf I} C \mySlice x_{i}
  \end{math}
  has the left adjoint
  \begin{equation*}
    \myLim_{i \myElemOf I} C \mySlice x_{i} \mySub C^{I} \mySlice x \myXMorphism{\myColim} C \mySlice \myColim_{i \myElemOf I} x_{i}.
  \end{equation*}
  The unit is an equivalence by the ``only if'' direction of the descent property.
  The counit is an equivalence by the ``if'' direction of the descent property.
\end{myProof}

\begin{myProposition}[{cf. \myTextCite[C3.3.1]{johnstone2002sketches2}}]
  \label{prop-prsv-dp-equiv-idx-la}
  Let \(F \myElemOf C \myMorphism D\) be a morphism of logoi.
  Then \(F\) preserves dependent products if and only if \(F\) has a \(C\)-indexed left adjoint.
\end{myProposition}
\begin{myProof}
  If \(F\) has a left adjoint \(G\), then \(F\) preserves dependent products if and only if \(G\) is \(C\)-indexed by \cref{prop-left-right-bc}.
  This implies the ``if'' direction.
  For the ``only if'' direction, it is enough to check \(F\) has a left adjoint.
  Since \(F\) preserves small colimits and finite limits, by \cref{prop-pres-right-aft}, it suffices to show that \(F\) preserves small products.
  Let \(I\) be a small type.
  The diagonal functor \(C \myMorphism C^{I}\) factors as the pullback functor \(C \myMorphism C \mySlice \myCoprod_{\myBlank \myElemOf I} \myTerminal\) followed by the equivalence \(C \mySlice \myCoprod_{\myBlank \myElemOf I} \myTerminal \myEquiv C^{I}\) (\cref{prop-logos-self-index-sheaf}) defined by the pullbacks along the inclusions \(\myTerminal \myMorphism \myCoprod_{\myBlank \myElemOf I} \myTerminal\) for all \(i \myElemOf I\).
  Thus, the product of a family \(x \myElemOf I \myMorphism C\) is constructed as the dependent product of \(\myCoprod_{i \myElemOf I} x_{i} \myElemOf C \mySlice \myCoprod_{\myBlank \myElemOf I} \myTerminal\) along \(\myCoprod_{\myBlank \myElemOf I} \myTerminal \myMorphism \myTerminal\).
  By assumption, \(F\) preserves this construction.
\end{myProof}

\begin{myRemark}
  A morphism of topoi corresponding to a morphism of logoi that has an indexed left adjoint is called a locally contractible morphism \myCite[Section 4.2.7]{anel2021topo-logie}, which is an \(\infty\)-analogue of a locally connected morphism in \(1\)-topos theory \myCite[C3.3.1]{johnstone2002sketches2}.
\end{myRemark}

\section{Univalence in logoi}
\label{sec:univalence-logoi}

We recall the notion of a univalent family in a logos.
Note that univalence is defined in more general context \myCite{gepner2017univalence,rasekh2021univalence,nguyen2025type}.

\begin{myDefinition}
  Let \(C\) be a logos.
  Let \(\myFamily(C) \mySub \myArrowCat{C}\) denote the wide subcategory spanned by the pullback squares.
  We call an object in \(\myFamily(C)\) a \myDefine{family in \(C\)}.
  For a family \(u\) in \(C\), its codomain is referred to as \(\myFamBase(u)\).
  A family is \myDefine{univalent} if it is a \((-1)\)-truncated object in \(\myFamily(C)\).
  Let \(\myUnivFam(C) \mySub \myFamily(C)\) denote the full subcategory spanned by the univalent families.
\end{myDefinition}

We show that the inclusion \(\myUnivFam(C) \myMorphism \myFamily(C)\) has a left adjoint (\cref{prop-univalent-completion}).
This can be proved by using object classifiers \myCite[Theorem 6.1.6.8]{lurie2009higher}, but we give an elementary construction.
The idea comes from the construction of \((-1)\)-truncation by iterated join in homotopy type theory \myCite{rijke2017join-arxiv}.

\begin{myLemma}
  \label{prop-fam-slice}
  Let \(C\) be a logos and \(u\) a family in \(C\).
  Then the functor
  \begin{equation*}
    \myFamBase \mySlice u \myElemOf \myFamily(C) \mySlice u \myMorphism C \mySlice \myFamBase(u)
  \end{equation*}
  is an equivalence.
\end{myLemma}
\begin{myProof}
  The inverse sends \(f \myElemOf C \mySlice \myFamBase(u)\) to the pullback of \(u\) along \(f\).
\end{myProof}

\begin{myLemma}
  \label{prop-fam-colimit}
  Let \(C\) be a logos.
  Then \(\myFamily(C) \mySub \myArrowCat{C}\) is closed under small colimits and pullbacks.
\end{myLemma}
\begin{myProof}
  We apply \cref{prop-closed-under-colimit-criterion-1} to see that \(\myFamily(C) \mySub \myArrowCat{C}\) is closed under small colimits.
  \Cref{cond-closed-under-colimit-criterion-1-proj} follows from the descent property.
  For every \(u \myElemOf \myFamily(C)\), the inclusion \(\myFamily(C) \mySlice u \myMorphism \myArrowCat{C} \mySlice u\) is equivalent by \cref{prop-fam-slice} to the functor \(C \mySlice \myFamBase(u) \myMorphism \myArrowCat{C} \mySlice u\) obtained by the pullback along \(u\), which preserves small colimits by \cref{prop-logos-colimit-pb-stable}.
  \Cref{cond-closed-under-colimit-criterion-1-slice} thus holds.
  Since pullbacks are binary products in slices and the pullback functor \(C \mySlice \myFamBase(u) \myMorphism \myArrowCat{C} \mySlice u\) preserves arbitrary limits, \(\myFamily(C) \mySub \myArrowCat{C}\) is closed under pullbacks.
\end{myProof}

\begin{myLemma}
  \label{prop-fam-bin-prod}
  Let \(C\) be a logos.
  Then \(\myFamily(C)\) has binary products that preserve small colimits on each variable.
\end{myLemma}
\begin{myProof}
  Let \(u_{1}\) and \(u_{2}\) be families in \(C\).
  A span \(u_{1} \myMorphismOp v \myMorphism u_{2}\) in \(\myFamily(C)\) corresponds to a span \(\myFamBase(u_{1})\myXMorphismOp{f_{1}} x \myXMorphism{f_{2}} \myFamBase(u_{2})\) in \(C\) equipped with an equivalence \(e \myElemOf \myPBFun{f_{1}}(u_{1}) \myEquiv \myPBFun{f_{2}}(u_{2})\) in \(C \mySlice x\).
  Let \(P\) be the presheaf on \(C \mySlice (\myFamBase(u_{1}) \myBinProd \myFamBase(u_{2}))\) that sends \((x \myComma (f_{1} \myComma f_{2}))\) to the type of equivalences \(\myPBFun{f_{1}}(u_{1}) \myEquiv \myPBFun{f_{2}}(u_{2})\) in \(C \mySlice x\).
  By \cref{prop-logos-self-index-sheaf}, \(P\) preserves small limits and thus is representable by \cref{prop-pres-left-repr}.
  The universal element for \(P\) then defines the product of \(u_{1}\) and \(u_{2}\).
  We show that the functor \((u_{1} \myBinProd \myHole) \myElemOf \myFamily(C) \myMorphism \myFamily(C)\) preserves small colimits.
  Since the forgetful functor from a slice creates colimits, it suffices to see that \((u_{1} \myBinProd \myHole) \mySlice v \myElemOf \myFamily(C) \mySlice v \myMorphism \myFamily(C) \mySlice (u_{1} \myBinProd v)\) preserves small colimits for every \(v \myElemOf \myFamily(C)\).
  We show that the corresponding functor \(C \mySlice \myFamBase(v) \myMorphism C \mySlice \myFamBase(u_{1} \myBinProd v)\) by \cref{prop-fam-slice} is the pullback functor along the projection \(\myFamBase(u_{1} \myBinProd v) \myMorphism \myFamBase(v)\), which preserves small colimits by \cref{prop-logos-colimit-pb-stable}.
  That is, for every morphism \(f \myElemOf u_{2} \myMorphism v\) in \(\myFamily(C)\), the square
  \begin{equation*}
    \begin{tikzcd}
      \myFamBase(u_{1} \myBinProd u_{2})
      \arrow[r, "\myFamBase(u_{1} \myBinProd f)"]
      \arrow[d] &
      [4ex]
      \myFamBase(u_{1} \myBinProd v)
      \arrow[d] \\
      \myFamBase(u_{2})
      \arrow[r, "\myFamBase(f)"'] &
      \myFamBase(v)
    \end{tikzcd}
  \end{equation*}
  is a pullback in \(C\).
  By definition, a morphism \(x \myMorphism \myFamBase(u_{1} \myBinProd u_{2})\) in \(C\) corresponds to a tuple \((g_{1} \myComma g_{2} \myComma e)\) of morphisms \(g_{1} \myElemOf x \myMorphism \myFamBase(u_{1})\) and \(g_{2} \myElemOf x \myMorphism \myFamBase(u_{2})\) in \(C\) and an equivalence \(e \myElemOf \myPBFun{g_{1}}(u_{1}) \myEquiv \myPBFun{g_{2}}(u_{2})\) in \(C \mySlice x\).
  Since \(u_{2} \myEquiv \myPBFun{\myFamBase(f)}(v)\), the pullback of \(\myFamBase(u_{1} \myBinProd v)\) along \(\myFamBase(f)\) has the same universal property.
\end{myProof}

\begin{myRemark}
  From the proof of \cref{prop-fam-bin-prod}, the inclusion \(\myFamily(C) \myMorphism \myArrowCat{C}\) does not preserve binary products.
\end{myRemark}

\begin{myCorollary}
  \label{prop-univ-bin-prod-prsv-omega-colim}
  Let \(C\) be a logos and \(u \myElemOf \myFinOrd \myMorphism \myFamily(C)\) be a diagram, where \(\myFinOrd\) is the poset of finite ordinals.
  Then the canonical morphism
  \begin{math}
    \myColim_{n \myElemOf \myFinOrd} (u_{n} \myBinProd u_{n}) \myMorphism (\myColim_{n \myElemOf \myFinOrd} u_{n}) \myBinProd (\myColim_{m \myElemOf \myFinOrd} u_{m})
  \end{math}
  is an equivalence.
  Consequently, \(\myUnivFam(C) \mySub \myFamily(C)\) is closed under colimits indexed by \(\myFinOrd\).
\end{myCorollary}
\begin{myProof}
  By
  \begin{align*}
    & (\myColim_{n \myElemOf \myFinOrd} u_{n}) \myBinProd (\myColim_{m \myElemOf \myFinOrd} u_{m})
    \\ \myEquiv{}
    & \myColim_{n \myComma m \myElemOf \myFinOrd} (u_{n} \myBinProd u_{m}) \tag{\cref{prop-fam-bin-prod}}
    \\ \myEquiv{}
    & \myColim_{n \myElemOf \myFinOrd} (u_{n} \myBinProd u_{n}) \tag{\(\myFinOrd \myMorphism \myFinOrd \myBinProd \myFinOrd\) is final}.
  \end{align*}
  The last assertion is because \(u \myElemOf \myFamily(C)\) is univalent if and only if the diagonal morphism \(u \myMorphism u \myBinProd u\) is an equivalence.
\end{myProof}

\begin{myProposition}
  \label{prop-univalent-completion}
  Let \(C\) be a logos.
  Then the inclusion \(\myUnivFam(C) \myMorphism \myFamily(C)\) has a left adjoint.
\end{myProposition}
\begin{myProof}
  We use a version of join construction \myCite{rijke2017join-arxiv}.
  Remember that \(\myFamily(C)\) has small colimits (\cref{prop-fam-colimit}) and binary products (\cref{prop-fam-bin-prod}).
  Let \(u\) be a family in \(C\).
  We construct a chain \(v \myElemOf \myFinOrd \myMorphism \myFamily(C)\) as follows.
  We set \(v_{0} \myDefEq u\).
  We define the successor step \(i_{n} \myElemOf v_{n} \myMorphism v_{n + 1}\) by the join of two copies of \(v_{n}\), that is, the pushout
  \begin{equation*}
    \begin{tikzcd}
      v_{n} \myBinProd v_{n}
      \arrow[r, "\myLimProj{2}"]
      \arrow[d, "\myLimProj{1}"']
      \arrow[dr, myPOMark] &
      v_{n}
      \arrow[d, "i'_{n}"] \\
      v_{n}
      \arrow[r, "i_{n}"'] &
      v_{n + 1}.
    \end{tikzcd}
  \end{equation*}
  Note that \(i_{n}\) and \(i'_{n}\) are identical to the composite \(v_{n} \myMorphism v_{n} \myBinProd v_{n} \myMorphism v_{n + 1}\).
  Let \(v_{\infty}\) be the colimit of \(v\).
  By construction, the morphism \(v_{n} \myBinProd v_{n} \myMorphism v_{n + 1}\) is the following diagonal filler.
  \begin{equation*}
    \begin{tikzcd}
      v_{n}
      \arrow[r, "i_{n}"]
      \arrow[d] &
      v_{n + 1}
      \arrow[d] \\
      v_{n} \myBinProd v_{n}
      \arrow[r, "i_{n} \myBinProd i_{n}"']
      \arrow[ur, dotted] &
      v_{n + 1} \myBinProd v_{n + 1}
    \end{tikzcd}
  \end{equation*}
  Then \(v_{\infty} \myEquiv v_{\infty} \myBinProd v_{\infty}\) by \cref{prop-univ-bin-prod-prsv-omega-colim}, and thus \(v_{\infty}\) is \((-1)\)-truncated.
  To see that \(v_{\infty}\) is the least univalent family with a morphism \(u \myMorphism v_{\infty}\), let \(v'\) be a univalent family with a morphism \(u \myMorphism v'\).
  We can extend the morphism to \(v'\) along each step \(i_{n} \myElemOf v_{n} \myMorphism v_{n + 1}\) because the square
  \begin{equation*}
    \begin{tikzcd}
      v_{n} \myBinProd v_{n}
      \arrow[r, "\myLimProj{2}"]
      \arrow[d, "\myLimProj{1}"'] &
      v_{n}
      \arrow[d] \\
      v_{n}
      \arrow[r] &
      v'
    \end{tikzcd}
  \end{equation*}
  trivially commutes as \(v'\) is \((-1)\)-truncated.
  Taking the colimit, we have \(v_{\infty} \myLe v'\).
\end{myProof}

We call the left adjoint of the inclusion \(\myUnivFam(C) \myMorphism \myFamily(C)\) the \myDefine{univalent completion functor}.

\begin{myRemark}
  The term ``univalent completion'' is also used by \myTextCite{vandenberg2018completion} to mean embedding a fibration into a univalent one via a homotopy pullback in the Kan-Quillen model structure on the category of simplicial sets.
  Our univalent completion in the logos of types gives univalent completion in their sense.
\end{myRemark}

\begin{myRemark}
  The construction of univalent completion given in the proof of \cref{prop-univalent-completion} also shows that univalent completion in a logos is stable under pullback.
\end{myRemark}

Univalence is preserved by a functor preserving dependent products \myCite[cf.][Corollary 5.11]{nguyen2025type}.

\begin{myProposition}
  \label{prop-prsv-dep-prod-prsv-univ}
  Let \(F \myElemOf C \myMorphism D\) be a morphism of logoi.
  If \(F\) preserves dependent products, then \(F\) takes univalent families in \(C\) to univalent families in \(D\).
\end{myProposition}
\begin{myProof}
  Let \(u\) be a univalent family in \(C\) and \(v\) a family in \(D\).
  By \cref{prop-prsv-dp-equiv-idx-la}, the equivalence \(\myHom_{\myArrowCat{D}}(v \myComma F(u)) \myEquiv \myHom_{\myArrowCat{C}}(G(v) \myComma u)\) restricts to a \((-1)\)-truncated map \(\myHom_{\myFamily(D)}(v \myComma F(u)) \myMorphism \myHom_{\myFamily(C)}(G(v) \myComma u)\), where \(G\) is the left adjoint of \(F\).
  Since \(\myHom_{\myFamily(C)}(G(v) \myComma u)\) is \((-1)\)-truncated by assumption, \(\myHom_{\myFamily(D)}(v \myComma F(u))\) is also \((-1)\)-truncated.
  Hence, \(F(u)\) is univalent.
\end{myProof}

\section{Étale morphisms of topoi}
\label{sec:etale-morph-topoi}

\begin{myProposition}[{\myTextCite[Proposition 6.3.5.1]{lurie2009higher}}]
  Let \(C\) be a logos and \(x \myElemOf C\).
  Then \(C \mySlice x\) is a logos, and the pullback functor \(\myPBFun{x} \myElemOf C \myMorphism C \mySlice x\) is a morphism of logoi.
\end{myProposition}
\begin{myProof}
  \(\mySumFun{x} \myElemOf C \mySlice x \myMorphism C\) creates small colimits and pullbacks.
  The identity on \(x\) is the terminal object in \(C \mySlice x\).
  Thus, \(C \mySlice x\) has small colimits and finite limits, and the descent property is inherited from \(C\).
  By \cref{prop-pres-slice}, \(C \mySlice x\) is presentable.
  \(\myPBFun{x}\) is a morphism of logoi by \cref{prop-logos-colimit-pb-stable}.
\end{myProof}

\begin{myDefinition}
  We say a morphism \(F \myElemOf C \myMorphism D\) of logoi is a \myDefine{base change morphism} if there exist an object \(x \myElemOf C\) and an equivalence \(F \myEquiv \myPBFun{x}\) in \(C \mySlice \myLogos\).
  Since \(C \myEquiv C \mySlice \myTerminal\) and \((C \mySlice x) \mySlice y \myEquiv C \mySlice \mySumFun{x}(y)\) for \(y \myElemOf C \mySlice x\), base change morphisms of logoi are closed under composition.
  Let \(\myLogosBC \mySub \myLogos\) denote the wide subcategory spanned by the base change morphisms.
  Let \(\myToposEt \myDefEq \myOpCat{\myLogosBC}\).
  A morphism in \(\myToposEt\) is called an \myDefine{étale morphism} of topoi.
\end{myDefinition}

The goal of this section is to prove \cref{prop-etale-colimit}, which asserts that \(\myToposEt \mySub \myTopos\) is closed under small colimits.
Our proof is based on \cref{prop-bc-char-1}, a characterization of base change morphisms of logoi in relation to univalent families.

\begin{myDefinition}
  We say a morphism \(F \myElemOf C \myMorphism D\) of logoi \myDefine{provides enough families} if, for every family \(v\) in \(D\), there exist a family \(u\) in \(C\) and a morphism \(v \myMorphism F(u)\) of families in \(D\).
\end{myDefinition}

\begin{myDefinition}
  Let \(F \myElemOf C \myMorphism D\) be a morphism of logoi that preserves dependent products.
  We say \(F\) \myDefine{provides enough univalent families} if, for every univalent family \(v\) in \(D\), there exists a univalent family \(u\) in \(C\) such that \(v \myLe F(u)\).
  Note that \(F(u)\) is univalent by \cref{prop-prsv-dep-prod-prsv-univ}.
\end{myDefinition}

\begin{myDefinition}
  We say a morphism \(F \myElemOf C \myMorphism D\) of logoi is \myDefine{essentially surjective on objects} if the closure of the image of \(F\) under small colimits and finite limits is \(D\).
\end{myDefinition}

We prepare some lemmas for \cref{prop-bc-char-1}.

\begin{myLemma}
  \label{prop-etale-transpose-pb}
  Let \(C\) be a logos and \(x \myElemOf C\) an object.
  Then, a commutative square in \(C \mySlice x\) of the form
  \begin{equation}
    \label[square]{eq-square-to-pb}
    \begin{tikzcd}
      z'
      \arrow[r]
      \arrow[d, "q"'] &
      \myPBFun{x}(y')
      \arrow[d, "\myPBFun{x}(p)"] \\
      z
      \arrow[r] &
      \myPBFun{x}(y)
    \end{tikzcd}
  \end{equation}
  is a pullback if and only if its transpose
  \begin{equation}
    \label[square]{eq-square-from-sum}
    \begin{tikzcd}
      \mySumFun{x}(z')
      \arrow[r]
      \arrow[d, "\mySumFun{x}(q)"'] &
      y'
      \arrow[d, "p"] \\
      \mySumFun{x}(z)
      \arrow[r] &
      y
    \end{tikzcd}
  \end{equation}
  is a pullback.
\end{myLemma}
\begin{myProof}
  Consider the following diagram in \(C\).
  \begin{equation*}
    \begin{tikzcd}
      \mySumFun{x}(z')
      \arrow[r]
      \arrow[d, "\mySumFun{x}(q)"'] &
      x \myBinProd y'
      \arrow[r]
      \arrow[d]
      \arrow[dr, myPBMark] &
      y'
      \arrow[d, "p"] \\
      \mySumFun{x}(z)
      \arrow[r]
      \arrow[dr] &
      x \myBinProd y
      \arrow[r]
      \arrow[d]
      \arrow[dr, myPBMark] &
      y
      \arrow[d] \\
      & x
      \arrow[r] &
      \myTerminal.
    \end{tikzcd}
  \end{equation*}
  The composite of the upper left and upper right squares is \cref{eq-square-from-sum}.
  \Cref{eq-square-to-pb} is represented by the upper left square.
  Since the upper right square is a pullback, the claim follows from the right cancellation property of pullbacks.
\end{myProof}

\begin{myCorollary}
  \label{prop-base-change-to-idx-la}
  Every base change morphism \(F \myElemOf C \myMorphism D\) of logoi has a \(C\)-indexed left adjoint and thus preserves dependent products by \cref{prop-prsv-dp-equiv-idx-la}.
\end{myCorollary}
\begin{myProof}
  By \cref{prop-etale-transpose-pb}.
\end{myProof}

\begin{myLemma}
  \label{prop-idx-la-factor-bc-ff}
  Let \(F \myElemOf C \myMorphism D\) be a morphism of logoi that has a \(C\)-indexed left adjoint \(G\).
  Then the morphism of logoi \(H \myElemOf C \mySlice G(\myTerminal) \myMorphism D\) defined by the composite
  \begin{equation*}
    C \mySlice G(\myTerminal) \myXMorphism{F \mySlice G(\myTerminal)} D \mySlice F(G(\myTerminal)) \myXMorphism{\myPBFun{\myAdjUnit}} D,
  \end{equation*}
  where \(\myAdjUnit\) is the unit of \(G \myAdjRel F\), is fully faithful and satisfies \(H \myComp \myPBFun{G(\myTerminal)} \myEquiv F\).
\end{myLemma}
\begin{myProof}
  By construction, \(H(\myPBFun{G(\myTerminal)}(x))\) for \(x \myElemOf C\) fits into the pullback
  \begin{equation*}
    \begin{tikzcd}
      H(\myPBFun{G(\myTerminal)}(x))
      \arrow[r]
      \arrow[d]
      \arrow[dr, myPBMark] &
      F(x) \myBinProd F(G(\myTerminal))
      \arrow[d] \\
      \myTerminal
      \arrow[r, "\myAdjUnit"'] &
      F(G(\myTerminal))
    \end{tikzcd}
  \end{equation*}
  and thus is equivalent to \(F(x)\).
  Observe that \(H\) has the left adjoint \(G \mySlice \myTerminal \myElemOf D \myMorphism C \mySlice G(\myTerminal)\).
  Its counit at \(x \myElemOf C \mySlice G(\myTerminal)\) is the transpose of the pullback
  \begin{equation*}
    \begin{tikzcd}
      H(x)
      \arrow[r]
      \arrow[d]
      \arrow[dr, myPBMark] &
      F(\mySumFun{G(\myTerminal)}(x))
      \arrow[d] \\
      \myTerminal
      \arrow[r, "\myAdjUnit"'] &
      F(G(\myTerminal))
    \end{tikzcd}
  \end{equation*}
  and thus an equivalence.
  Hence, \(H\) is fully faithful.
\end{myProof}

\begin{myRemark}
  A morphism of topoi corresponding to a fully faithful morphism of logoi is called contractible \myCite[Section 4.2.7]{anel2021topo-logie}, which is an \(\infty\)-analog of a connected morphism in \(1\)-topos theory \myCite[C1.5.7]{johnstone2002sketches2}.
  \Cref{prop-idx-la-factor-bc-ff} thus gives a factorization of a locally contractible morphism of topoi into a contractible morphism followed by an étale morphism \myCite[cf.][C3.3.5]{johnstone2002sketches2}.
\end{myRemark}

\begin{myProposition}
  \label{prop-bc-char-1}
  Let \(F \myElemOf C \myMorphism D\) be a morphism of logoi that has a \(C\)-indexed left adjoint \(G\) (or, equivalently by \cref{prop-prsv-dp-equiv-idx-la}, preserves dependent products).
  Then the following are equivalent.
  \begin{enumerate}
    \NewCommandCopy{\myLabelOrig}{\label}
    \RenewDocumentCommand{\label}{m}{\myLabelOrig[cond]{#1}}
  \item \label{cond-bc-char-1-bc} \(F\) is a base change morphism.
  \item \label{cond-bc-char-1-unit-cart} The unit of \(G \myAdjRel F\) is cartesian.
  \item \label{cond-bc-char-1-enough-fam} \(F\) provides enough families.
  \item \label{cond-bc-char-1-enough-univ-fam} \(F\) provides enough univalent families.
  \item \label{cond-bc-char-1-eso} \(F\) is essentially surjective on objects.
  \end{enumerate}
  Moreover, if these conditions are satisfied, then \(G \mySlice \myTerminal \myElemOf D \myMorphism C \mySlice G(\myTerminal)\) gives an equivalence \(F \myEquiv \myPBFun{G(\myTerminal)}\) in \(C \mySlice \myLogos\).
\end{myProposition}
\begin{myProof}
  \Cref{cond-bc-char-1-bc} implies \cref{cond-bc-char-1-unit-cart} by \cref{prop-etale-transpose-pb} as the unit \(y \myMorphism F(G(y))\) is the transpose of the identity.
  The implications
  \begin{math}
    \labelcref{cond-bc-char-1-unit-cart} \implies \labelcref{cond-bc-char-1-enough-fam}
  \end{math}
  and
  \begin{math}
    \labelcref{cond-bc-char-1-enough-fam} \implies \labelcref{cond-bc-char-1-eso}
  \end{math}
  are trivial.
  By \cref{prop-idx-la-factor-bc-ff}, \(F\) factors into a base change morphism \(K\) followed by a fully faithful morphism \(H\).
  If \(F\) is essentially surjective on objects, then so is \(H\).
  Observe that a morphism of logoi that is both fully faithful and essentially surjective on objects is an equivalence.
  Therefore, \cref{cond-bc-char-1-eso} implies \cref{cond-bc-char-1-bc}.
  \Cref{cond-bc-char-1-enough-fam,cond-bc-char-1-enough-univ-fam} are equivalent by univalent completion (\cref{prop-univalent-completion}).
  The last assertion is because \(G \mySlice \myTerminal\) is the left adjoint of \(H\) as in \cref{prop-idx-la-factor-bc-ff}.
\end{myProof}

\begin{myCorollary}[{\myTextCite[Remark 6.3.5.7]{lurie2009higher}}]
  \label{prop-logos-embed-coslice-bc}
  Let \(C\) be a logos.
  Then the functor
  \begin{math}
    x \myMapsTo \myPBFun{x} \myElemOf \myOpCat{C} \myMorphism C \mySlice \myLogosBC
  \end{math}
  is an equivalence.
\end{myCorollary}
\begin{myProof}
  Define a functor \(L \myElemOf C \mySlice \myLogosBC \myMorphism \myOpCat{C}\) by \(L(F) \myDefEq G(\myTerminal)\),  where \(G\) is the \(C\)-indexed left adjoint of \(F \myElemOf C \myMorphism D\) (\cref{prop-base-change-to-idx-la}).
  By \cref{prop-bc-char-1}, \(F \myEquiv \myPBFun{L(F)}\) for \(F \myElemOf C \mySlice \myLogosBC\).
  We also have \(x \myEquiv L(\myPBFun{x})\) for \(x \myElemOf \myOpCat{C}\) by construction.
\end{myProof}

We turn to the proof of \cref{prop-etale-colimit}.

\begin{myLemma}[{\myTextCite[Proposition 6.3.5.14]{lurie2009higher}}]
  \label{prop-coslice-limit-base-change}
  Let \(C\) be a logos.
  Then \(C \mySlice \myLogosBC \mySub C \mySlice \myLogos\) is closed under small limits.
\end{myLemma}
\begin{myProof}
  By \cref{prop-logos-embed-coslice-bc}, it suffices to show that the functor
  \begin{math}
    x \myMapsTo \myPBFun{x} \myElemOf \myOpCat{C} \myMorphism C \mySlice \myLogos
  \end{math}
  preserves small limits, but this follows from \cref{prop-logos-self-index-sheaf}.
\end{myProof}

\begin{myLemma}
  \label{prop-prod-proj-base-change}
  Let \(C \myElemOf I \myMorphism \myLogos\) be a small family of logoi.
  Then the projection \(\myLimProj{i} \myElemOf \myProd_{j \myElemOf I} C_{j} \myMorphism C_{i}\) is a base change morphism for every \(i \myElemOf I\).
\end{myLemma}
\begin{myProof}
  We define an object \(x \myElemOf \myProd_{j \myElemOf I} C_{j}\) by \(x_{j} \myDefEq \myCoprod_{\myBlank \myElemOf j \myId i} \myTerminal\).
  Then
  \begin{equation*}
    (\myProd_{j \myElemOf I} C_{j}) \mySlice x \myEquiv \myProd_{j \myElemOf I} C_{j} \mySlice x_{j} \myEquiv \myProd_{j \myElemOf I} \myProd_{\myBlank \myElemOf j \myId i} C_{j} \myEquiv C_{i},
  \end{equation*}
  where the second equivalence follows from \cref{prop-logos-self-index-sheaf}.
  The composite
  \begin{math}
    \myProd_{j \myElemOf I} C_{j} \myXMorphism{\myPBFun{x}} (\myProd_{j \myElemOf I} C_{j}) \mySlice x \myEquiv C_{i}
  \end{math}
  coincides with the projection.
\end{myProof}

\begin{myLemma}
  \label{prop-lim-proj-base-change}
  Let \(C \myElemOf I \myMorphism \myLogosBC\) be a small diagram, and take the limit \(C_{-\infty}\) of \(C\) in \(\myLogos\).
  Then the projection \(\myLimProj{i} \myElemOf C_{-\infty} \myMorphism C_{i}\) is a base change morphism for every \(i \myElemOf I\).
\end{myLemma}
\begin{myProof}
  Let \(D \myDefEq \myProd_{j \myElemOf I} C_{j}\).
  By \cref{prop-prod-proj-base-change}, it suffices to show that the forgetful morphism \(\pi \myElemOf C_{-\infty} \myMorphism D\) is a base change morphism.
  We apply \cref{prop-bc-char-1}.
  By \cref{prop-logos-lcc-limit,prop-base-change-to-idx-la}, \(\pi\) preserves dependent products.
  We show that \(\pi\) provides enough univalent families.
  Observe that \(\myUnivFam(C_{-\infty})\) is the full subposet of \(\myUnivFam(D)\) spanned by those \(u\) such that \(C_{s}(u_{i}) \myId u_{j}\) for every \(s \myElemOf i \myMorphism j\) in \(I\).
  We thus show that, for every \(u \myElemOf \myUnivFam(D)\), there exists a \(u^{\infty} \myElemOf \myUnivFam(C_{-\infty})\) such that \(u \myLe u^{\infty}\).
  We find such a \(u^{\infty}\) by Adámek's fixed point theorem \myCite{adamek1974free}.

  For a morphism \(s \myElemOf i \myMorphism j\) in \(I\), let \(C_{s}^{!}\) denote the left adjoint of \(C_{s} \myElemOf C_{i} \myMorphism C_{j}\).
  By \cref{prop-fam-colimit}, \(\myFamily(D)\) has small coproducts. \(C_{s}^{!}\) preserves pullbacks as \(\mySumFun{x}\) does for any \(x \myElemOf C_{i}\).
  We then define a functor \(F \myElemOf \myFamily(D) \myMorphism \myFamily(D)\) by
  \begin{equation*}
    F(u)_{i} \myDefEq (\myCoprod_{(j \myComma s) \myElemOf I \mySlice i} C_{s}(u_{j})) \myBinCoprod (\myCoprod_{(j \myComma s) \myElemOf i \mySlice I} C_{s}^{!}(u_{j})).
  \end{equation*}
  Let \(G\) be the composite
  \begin{equation*}
    \myUnivFam(D) \mySub \myFamily(D) \myXMorphism{F} \myFamily(D) \myMorphism \myUnivFam(D),
  \end{equation*}
  where the last one is the univalent completion functor (\cref{prop-univalent-completion}).
  By construction, we have morphisms
  \begin{align}
    & \label{eq-cr-to-g} C_{s}(u_{i}) \myMorphism G(u)_{j} \\
    & \label{eq-cl-to-g} C_{s}^{!}(u_{j}) \myMorphism G(u)_{i}
  \end{align}
  in \(\myFamily(D)\) for every morphism \(s \myElemOf i \myMorphism j\) in \(I\).
  Note that the latter morphism corresponds to a morphism \(u_{j} \myMorphism C_{s}(G(u)_{i})\) in \(\myFamily(D)\) by \cref{prop-etale-transpose-pb}.
  As a special case of the former, we have \(u \myLe G(u)\).
  For \(u \myElemOf \myUnivFam(D)\), take the colimit \(u^{\infty}\) of the chain
  \begin{equation*}
    \begin{tikzcd}
      u \myLe G(u) \myLe G^{2}(u) \myLe \dots
    \end{tikzcd}
  \end{equation*}
  indexed by \(\myFinOrd\).
  By \cref{prop-univ-bin-prod-prsv-omega-colim} and the construction of \(F\), the functor \(G\) preserves colimits indexed by \(\myFinOrd\), and thus \(u^{\infty} \myId G(u^{\infty})\).
  For any morphism \(s \myElemOf i \myMorphism j\) in \(I\), we have \(C_{s}(u^{\infty}_{i}) \myLe u^{\infty}_{j}\) by \cref{eq-cr-to-g} and \(u^{\infty}_{j} \myLe C_{s}(u^{\infty}_{i})\) by \cref{eq-cl-to-g}, and thus \(C_{s}(u^{\infty}_{i}) \myId u^{\infty}_{j}\).
  Therefore, \(u^{\infty}\) belongs to \(\myUnivFam(C_{-\infty})\).
\end{myProof}

\begin{myTheorem}[{\myTextCite[Theorem 6.3.5.13]{lurie2009higher}}]
  \label{prop-etale-colimit}
  \(\myToposEt \mySub \myTopos\) is closed under small colimits.
\end{myTheorem}
\begin{myProof}
  We apply \cref{prop-closed-under-colimit-criterion-1}.
  \Cref{cond-closed-under-colimit-criterion-1-slice} follows from \cref{prop-coslice-limit-base-change}.
  \Cref{cond-closed-under-colimit-criterion-1-proj} follows from \cref{prop-lim-proj-base-change}.
\end{myProof}

\section*{Acknowledgements}

The author was partially supported by JST (JPMJMS2033).

\printbibliography

\end{document}